\documentclass[12pt]{amsart}
\usepackage{amsmath}
\usepackage{amscd}
\usepackage{graphics}
\usepackage{latexsym}

\oddsidemargin .in
\theoremstyle{plain}
\newtheorem{Thm}{Theorem}
\newtheorem{Lem}[Thm]{Lemma}

\newtheorem{Rm}[Thm]{Remark}

\theoremstyle{definition}

\theoremstyle{remark}

\def\C{{\mathbb C}}
\def\R{{\mathbb R}}
\def\P{{\mathbb P}}
\def\Z{{\mathbb Z}}

\def\s2x{\hbox{$S^2 \times S^2$}}

	\def\sqr#1#2{{\vcenter{\hrule height.#2pt
    		\hbox{\vrule width.#2pt height#1pt \kern#1pt
       		\vrule width.#2pt}\hrule height.#2pt}}}
	\def\square{\mathchoice\sqr67\sqr67\sqr{2.1}6\sqr{1.5}6}

\def\hckin#1{#1}
\def\hckout#1{}

\begin{document}

\title[]{ Transcendental submanifolds of  ${\mathbb R}{\mathbb P}^n$ }
\author{Selman Akbulut and Henry King} 
\thanks{First named author is partially supported byNSF grant DMS 9971440}
\address{Deptartment  of Mathematics, Michigan State University,  MI, 48824}
\email{akbulut@math.msu.edu }
\address {Department. of Mathematics\ University of Maryland,  Md, 20742 }
\email{hck@math.md.edu }
\date{\today}
\begin{abstract}
In this paper we give examples of closed smooth submanifolds of $\R\P^{n}$ which are isotopic to
nonsingular projective subvarieties of $\R\P^{n}$ but they can not be isotopic to the real parts of nonsingular complex projective subvarieties of $\C\P^{n}$.
\end{abstract}

\maketitle
\setcounter{section}{-1}
\section{Introduction}

 Let
  $j: \R^n \hookrightarrow \R\P^n$ be the canonical imbedding as a chart.  Real  algebraic sets in $\R^n$ are not in general real algebraic sets in $\R\P^n$. The Zariski closure of  the image of an algebraic set (under  $j$) usually has extra components at infinity. An algebraic subset of $\R^n $ which remains an algebraic set in $\R\P^n$ is called a {\it projectively closed algebraic set} (\cite{ak1}). 
Not every algebraic set is projectively closed. 
In general, isotoping a submanifold of  the projective space $\R\P^n$ to an algebraic subset is a much harder problem than the corresponding problem in the affine case ${\mathbb R}^n$.  In this paper we  produce a transcendental submanifold of $\R\P^n$ in the sense of  \cite{ak5}.  That is, we find a smooth submanifold of  $\R\P^n$ which is isotopic to a nonsingular projective algebraic subset,  but which can not be isotoped to the real part  of any  complex nonsingular algebraic subset of  $\C\P^n$.  This results generalizes the affine examples of \cite{ak5} to the projective case. We want to thank MSRI for giving us \hckin{the} opportunity to work together.

\section{Preliminaries}

\hckin{By a {\em closed} (sub)manifold we mean a compact (sub)manifold without boundary.}

Let $V$ be a \hckout{Zariski open }real (or complex) algebraic set defined over $\R$, and let
$R={\Z}_{2}$ (or  $R={\Z}$). Then we can 
define algebraic homology groups 
 $H_{*}^{A}(V;R)$
to be the subgroup of $H_{*}(V;R)$ generated by the
compact real (or complex) algebraic subsets of $V$    (cf. \cite {ak1}).
We define $H^{*}_{A}(V;R)$ to be the Poincar\'e duals of the groups 
$H_{*}^{A}(V;R)$ when defined. 
By the resolution of singularities  theorem (\cite{h}), $H_{*}^{A}(V;R)$ is also the subgroup
generated by the classes $g_{\ast }([S])$ where 
$g\colon S \to V$ is an entire rational function,
$S$ is a compact nonsingular real (or complex) algebraic set and $[S]$ is the fundamental class of $S$. 
Therefore even when $V$ is real, we can define 
$H_{*}^{A}(V;{\Z})$
to be the subgroup generated by $g_{\ast }([S])$ where $g:S\to V$ is an entire rational function 
from an oriented compact nonsingular real algebraic set and $[S]$ is the fundamental class of $S$.

Now let
$V\subset \R\P^{n}$ be a nonsingular projective
real algebraic set of dimension $v$, and  suppose its complexification $V_\C\subset \C\P^{n}$ is nonsingular.
 Let $j:V\hookrightarrow V_\C$ denote the inclusion.
Define $\bar {H}_{2k}^{A}(V_\C;{\Z})$ to be the subgroup of $H_{2k}^{A}(V_\C;\Z)$
generated by irreducible complex algebraic subsets defined over $\R$ 
with $k$-dimensional real parts. 
In other words it is generated by the complexification of $k$-dimensional real algebraic 
subsets of $V$ in $V_\C$.
Again by the resolution theorem, $\bar {H}_{2k}^{A}(V_\C;\Z)$ is generated by
the classes $g_{\ast }([L_\C])$, where 
$L_\C$ is an irreducible nonsingular complex projective algebraic set 
defined over $\R$ with nonempty real part
and $g\colon L_\C \to V_\C$ is a regular map defined over $\R$.
Let $\bar{H}_{A}^{2k}(V_\C;\Z)$ denote the Poincar\'e dual of $\bar {H}_{2v-2k}^{A}(V_\C;\Z)$.
Define

$$\bar{H}_{\C-alg}^{2k}(V;\Z)=j^{*}\bar{H}_{A}^{2k}(V_\C;\Z)$$
Let $\bar{H}_{\C-alg}^{2k}(V;\Z_2)$ to be the mod $2$ reduction of $\bar{H}_{\C-alg}^{2k}(V;\Z)$ (under the obvious coefficient homomorphism $\Z \to {\Z}_{2}$). 
Define the natural subgroup of $H_{A}^{2k}(V;\Z_{2})$
$$ H_{A}^{k}(V;\Z_{2})^{2}=\{\;\alpha ^{2} \; |\; \alpha \in H_{A}^{k}(V;\Z_{2})\; \}$$ 
  Recall that Theorem A (b)  of \cite{ak5}  relates these groups to each other:

\vspace{.1in}

\begin{Thm}   For all $k$ the following hold: 
$\bar{H}_{\C-alg}^{2k}(V;\Z_{2})=H_{A}^{k}(V;\Z_{2})^{2}$.  
\end{Thm}

\vspace{.15in}


Let $M\subset V$ be a closed smooth submanifold of a nonsingular algebraic set $V$. The problem of whether $M$
is isotopic to a nonsingular algebraic subset  of $V$ is an old one. If we allow stabilization (replacing $V$ by $V\times {\R}^{k}$ for sufficiently large $k$) the problem becomes  solvable (\cite{n} and \cite{t} for the affine case, \cite{k} for the projective case, and \cite{ak2} for the general case where there are obstructions). 

If we don't allow stabilization the problem becomes much harder, in which the complexification of $V$ begins to play an important role. In \cite {ak3} it was shown that every closed smooth
$M\subset \R^{n}$ can be isotoped to a nonsingular (topological) component of an algebraic subset of ${\Bbb R}^n$. More generally
in \cite{ak4}  it was shown that any immersed submanifold of $\R^n$ can be isotoped 
to a nonsingular algebraic subset of ${\Bbb R}^n$ 
if and only if $M$ is cobordant through immersions to an algebraic subset of ${\Bbb R}^n$. 
Here we need a very special case of these theorems, 
the one which allows us to isotop some submanifolds of $\R\P^{n}$ 
(the examples in the next section)  to projectively closed algebraic subsets. 
Obviously the following lemmas hold in  more general contexts, 
but to make the examples of  the main theorem in the next section transparent, we chose to state them in this special form, which is enough to prove the theorem.

\begin{Lem}
Every closed codimension one submanifold of $\R^k$ can be \hckin{$C^\infty$} approximated by a nonsingular projectively closed algebraic subset.
\end{Lem}

\proof
This for example is proven in \cite{k} (also Theorem 2.8.2 of  \cite{ak1}),  but can also be seen from Seifert's original proof  \cite{s}  by noting that the highest degree terms of the polynomial he constructs are a constant times $|x|^{2n}$ (clearly the zeros of such polynomials are projectively closed algebraic sets).
\endproof

\begin{Lem}
Let  $M^{m}\subset Y^{m+1}$ be closed smooth manifolds with  $M$ separating $Y$. Let $f\colon Y\to \R^k$ be a smooth imbedding.
Then $f(M)$ is isotopic to a nonsingular projectively closed subvariety of $\R^k$. In particular,  
by viewing $f(M)$ as a submanifold of $\R\P^k$ via the natural inclusion $\R^k\subset \R\P^k$, $f(M)$ is isotopic to a nonsingular projective algebraic subvariety of $\R\P^k$. 
\end{Lem}

\proof  By \cite{ak3} we can isotop $f(Y)$ to a nonsingular  topological  component $Y'$ of a real algebraic subvariety of $\R^k$. 
Let  $W\subset Y$ be one of the codimension zero components of $Y-M$, 
and let $Q^{k-1}\subset \R^{k}$ be the boundary of a small tubular neighborhood of $f(W)$ in $\R^k$. By Lemma 2 above we can  isotop $Q$ to a projectively closed nonsingular algebraic subvariety $Q'$ of  $\R^k$ \hckin{which is $C^1$ close to $Q$}. Then $V:=Y'\cap Q'$ gives the desired variety. 
\hckin{$V$ is isotopic to $f(M)$ by transversality.}
$V$ is projectively closed since $Q'$ is projectively closed.
\endproof

 \begin{figure}[ht]
 \begin{center}
    \includegraphics{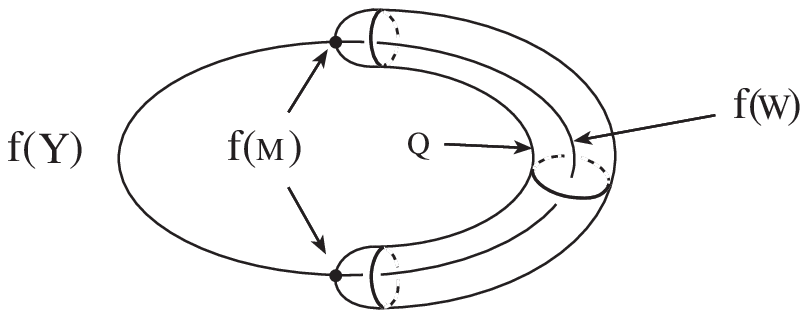}
  \caption{} \label{}
   \end{center}
 \end{figure}

\section{Transcendental submanifolds}

In \cite{ak5} we constructed smooth submanifolds $M\subset \R^{n}$ which are isotopic to nonsingular algebraic subsets of $\R^{n}$,  but not  isotopic to the real parts of nonsingular complex algebraic subsets of $\C\P^{n}$. So, this means that  either every algebraic model of $M$ in $ \R^{n}$ develops extra components at infinity when  its  Zariski closure is taken in $\R\P  ^n$ (i.e. it is not projectively closed), or   $M$ admits  nonsingular algebraic models in $\R\P ^{n}$  but the complexifications all such models necessarily contain singular points in $\C\P ^{n}$. The following theorem eliminates the first possibility, hence it gives genuine topological obstructions to moving smooth submanifolds of $\R\P^n$ to nonsingular algebraic sets in the strong sense (i.e. no singularities in complexification).

\begin{Thm} 
There are closed smooth submanifolds $M\subset \R\P^n$ which can be approximated (via a small  isotopy) by nonsingular subvarieties  of $ \R\P^n$,  but they can not be isotoped to
the real parts of nonsingular  complex algebraic subvarieties of $\C\P^n$ defined over $\R$. 
\end{Thm}

\begin{proof} By \cite{mm} for any $s$ there is an $m$ such that we have an imbedding 
$\R \P ^{m}\subset \R^{2m-s}$. Let $M= \R \P ^{m}\times S^{1} \subset \R^{2m-s}\times 
 \R ^{3}= \R^{2m+3-s}\subset \R\P^{n}$, where $n=2m+3-s$. Let  $Y=\R \P ^{m}\times S^{2}$ in $\R\P ^{n}$, 
so by Lemma 2 above $M$ can be isotoped to a nonsingular projectively closed algebraic subset $V^{v}$ of $\R\P ^{n}$
  where $v=m+1$.  
 
 We claim that when  $s\geq 3$ when $m$ is even, and $s\geq 5$  when $m=4k+1$, $V$  can not be  the real part of a nonsingular complex algebraic subset $V_{\C}$  of $\C\P^{n}$ (defined over $\R$).  Suppose such a $V_{\C}$ exists.
 By a version of the Lefschetz hyperplane theorem due to Larsen \cite {hr},
 for $i\leq 2v-n=s-1$, the  restriction induces an 
isomorphism
$$ H^{i}(\C\P^{n};\Z)\stackrel{\cong}{\to} H^{i}(V_\C;\Z)$$
  In particular when  $i\leq s-1$, the group $\bar{H}_{\C-alg}^{i}(V;\Z)$ lies in the the image of 
$H^{i}(\R\P^{n}; \Z)$ under restriction.
  But since $V$
lies in a chart $\R^{n}$ of $\R\P^{n}$, from the diagram
  
   $$\begin{array}{ccc}
  H^{i}(V_\C; \Z _{2}) & \stackrel{j^{*}}{\longrightarrow } 
    & H^{i}(V; \Z _{2}) \\
  \cong \uparrow &                   &  \uparrow  \\
  H^{i}(\C\P^{n}; \Z _{2}) 
  & \stackrel{j^{*}}{\longrightarrow } & H^{i}(\R\P ^{n}; \Z _{2})
\end{array}$$

\noindent we conclude that $\bar{H}_{\C-alg}^{i}(V;{\Z _{2}})=0$ for $0<i\leq s-1$.
On the other hand by \cite{ak1} the Stiefel-Whitney classes of $V$ are represented by real  algebraic subsets since $V$ is a nonsingular real algebraic set. Hence,  when $m$ is even $w_{1}(V)=\alpha  \times 1$ is algebraic, also when $m=4k+1$ then  $w_{2}(V)=\alpha ^{2} \times 1$ is algebraic, 
where $\alpha$ is the generator of $H^{1}(\R\P ^{m};\Z _{2})$.                                 
So when $m$ is even and $s\geq 3$,  or when $m=4k+1$ and $s\geq5$ we get a contradiction to Theorem 1 above, 
for example  when $m$ is even
$$0 \neq \alpha ^{2}\times 1= w_{1}^{2}(M)\in H_{  A}^{1 }(V;{\Z}_{2})^{2} = \bar{H}_{\C-alg}^{2}(V;{\Z}_{2})=0$$                  

We should point out that  the above mentioned theorem of \cite{mm} actually implies that there are imbeddings  $\R \P ^{m}\subset \R^{2m-s}$ for the pairs $(m,s)$ with  $m$ even and $s\geq 3$, or $m=4k+1$ and $s\geq 5$. This is what is used in the proof.             
\end{proof}                                                                         
                                                                                                                                                                                                             
\begin{Rm}
Recall that on any nonsingular real algebraic variety structure $V$ of a smooth manifold $M$,  the Stiefel-Whitney classes, and mod $2$ reductions of Pontryagin classes of $M$ are algebraic  in $V$  (\cite{ak1}, \cite{ak5}). So  the examples in the above  theorem generalize\hckout{s to} \hckin{in} many directions. For example we can take  $M^v$ to be any closed smooth manifold  which admits a separating imbedding into a closed manifold $Y^{v+1}$,  such that  $Y^{v+1}\subset \R^{2v-q}$ with $2i\leq q$, and $\alpha \in H^{i}(M:\Z_{2})$ such that  $\alpha^{2}\neq 0$ and  $\alpha$ lies in the subring of the cohomology group generated by Stiefel-Whitney and Pontryagin classes. 

\end{Rm}

\hckin{We thank the referees for helpful comments.}


\begin{thebibliography}{99999}


\bibitem[AK1]{ak1} S.~Akbulut and H.~King, The topology of real algebraic sets with isolated singularities, Ann of Math 113 (1981) 425-446

\bibitem[AK2]{ak2} S.~Akbulut and H.~King, The topology of real algebraic sets, M.S.R.I. book
series no.25 


\bibitem[AK3]{ak3} S.~Akbulut and H.~King, On approximating submanifolds by
algebraic sets and a solution to the Nash conjecture, Invent.\ Math.,107 (1992),87-98.

\bibitem[AK4]{ak4} S.~Akbulut and H.~King, Algebraicity of immersions, Topology, vol.31, no.4 (1992) 701-712.

\bibitem[AK5]{ak5} 
S. ~Akbulut and H.King,
{\em Transcendental submanifolds of ${\mathbb R}^n$}, Comm. Math. Helv 68 (1993), 308-318. 

\bibitem[Hr]{hr} R.~Hartshorne, Equivalence relations on algebraic cycles and subvarieties 
of small codimension, Proc. of Symposia in Pure Math. vol 29, 129-164 (1975)

\bibitem[H]{h} H.~Hironaka, Resolution of
singularities of an algebraic variety  over a field of characteristic zero,
Annals of Math.\ 79 (1964) 109-326.



\bibitem[K]{k} H. ~King,
	Approximating submanifolds of real projective space
	 by varieties,
            Topology  Vol. 15 (1976), pp. 81-85.

\bibitem[MM]{mm} M.~Mahowald and J.~Milgram, Embedding real projective spaces, Ann. of Math., 87 (1968), 411-422

\bibitem[N]{n} J. Nash, Real algebraic manifolds, 
 Annals of Math. 56 (1952), pp. 405-421.

\bibitem[S]{s} H. Seifert Algebraische approximation von mannigfaltigkeiten, Math.Zeit, 41 (1936), 1-17.


\bibitem[T]{t}  A.~Tognoli,   Su una congettura di Nash, Ann. Sc. Norm. Sup Pisa 27 (1973) pp. 167-185.
 \end{thebibliography}
\end{document}